\documentclass[11pt,notitlepage,twoside]{article}
\usepackage{latexsym}
\usepackage{amsmath}
\usepackage{amsfonts}
\usepackage{amssymb}
\renewcommand{\a }{\alpha }
\renewcommand{\b }{\beta }

\newcommand{\D }{\Delta }

\renewcommand{\l }{\lambda }

\renewcommand{\O }{\Omega }

\newcommand{\ov}{\overline}

\newcommand{\be}{\begin{equation}}
\newcommand{\ee}{\end{equation}}

\newcommand{\R}{\mathbb{R}}
 
\newcommand{\N}{\mathbb{N}}
\newtheorem{thm}{Theorem}[section]
\newtheorem{pro}[thm]{Proposition}
\newtheorem{lem}[thm]{Lemma}

\newtheorem{cor}[thm]{Corollary}

\numberwithin{equation}{section}

\author{{\bf\hfill\large Abdellaziz Harrabi\hfill}\vspace{1mm}\\
{\it\small Institut de Math{\'e}matiques Appliqu\'ees et d'Informatiques}\\
{\it\small  Kairouan, Tunisia}\\
{\it\small e-mail: Abdellaziz.Harrabi@ismai.rnu.tn}\vspace{1mm}\\}

\begin{document}

\date{ }
\title{Ambrosetti-Rabinowitz theorems revisited}
\maketitle

 {\footnotesize \noindent \ \ {\bf Abstract.} We introduce tow assumptions weaker than the
 classical Ambrosetti-Rabinowitz and the subcritical polynomial growth conditions
 to obtain the Palais-Smale Condition. Therefore, we improve the Ambrosetti-Rabinowitz existence theorems.
 Also, we prove some existence results under a weaken
 subcritical polynomial growth conditions where the nonlinearity does not satisfy the super-quadratic condition and even it does not have any limits at infinity.

}

   \section{ \bf Introduction and main results}
\mbox{}\quad This work is devoted to the following higher order elliptic equation under the
Navier boundary condition
 \begin{equation}
\label{a} \left\{\begin{array}{lll} (-\Delta)^m u=f(x,u)\ \ \   \mbox{in}\,\, \O,\\
\,\ \ u=\D u= . .=(\D)^{m-1}=0\ \ \ \mbox{on}\ \ \partial \O,
\end{array}
\right.
\end{equation}
or under the Dirichlet boundary condition
\begin{equation}
\label{aa} \left\{\begin{array}{lll} (-\Delta)^m u=f(x,u)\ \ \   \mbox{in}\,\, \O,\\
\,\ \ u=\nabla u= . .=\nabla^{m-1}u=0\ \ \ \mbox{on}\ \ \partial \O,
\end{array}
\right.
\end{equation}
where $\O$ is a bounded smooth domain in $\R^N$, $N\geq 2$.
Since $60$s, many authors have studied the existence or compactness of
solutions to \eqref{a} or \eqref{aa}. One classical approach is to use the variational method and the
 critical points theory (see for examples \cite{GGS, MS}). To make in work this method, we choose the following
 functional spaces: For equation \eqref{a}, let
$$ H_\vartheta^m(\Omega):=\left\{v\in H^m(\Omega);\;\Delta^jv=0\;\mbox{on}\;\partial\Omega,\;\mbox{for}\;
j<\frac{m}{2}\right\}.$$
For \eqref{aa}, let
$$ H_0^m(\Omega):=\left\{v\in
H^m(\Omega);\; \nabla^ju=0\;\mbox{on}\;\partial\Omega,\;\mbox{for}\;j=0, 1, . . ,m-1\right\}.$$
We denote by $E_m=H_\vartheta^m(\Omega)$ or $E_m=H_ 0^m(\Omega)$ and by $E'_m$ its dual space, we endow the Hilbert space $E_m$ with the standard
scalar product
\begin{equation}\label{1.7} (u,v )_m=\left\{\begin{array}{lll}&\displaystyle \int_\Omega(\Delta^l u)(\Delta^l v),\vspace*{2mm}&\quad\mbox{if } m=2l,\\
&\displaystyle \int_\Omega(\nabla\Delta^l u)(\nabla\Delta^l v),&\quad\mbox{if } m=2l+1.\end{array}\right. \end{equation}
We denote by $|\cdot|_{m}$ the corresponding norm. Let $f \in C(\overline\O\times\R)$. Under the following weak subcritical assumption :
 $$\mbox{There exist $C_0>0$, $s_0>0$ such that $|f(x,s)|\leq C_0|s|^{\frac{N+2m}{N-2m}}$,
$\forall\; |s|\geq s_0$ and $x\in \Omega$,} \leqno{(H)}$$
the energy functional

    $$ J(u)=\frac{1}{2}|u|_{m}^2-\int_\Omega
    F(x, u)$$
     is well defined and belongs to $C^1(E_m),$ here $F(x,s)= \int^s_0 f(x,t)dt$ for $(x,s) \in \O\times \R$.

   We say that $u\in E_m$ is a weak solution of problem \eqref{a}
 (resp.~\eqref{aa}) if
$$( u,v )_m= \displaystyle\int_\Omega
    {f}(u)v, \;\; \forall\; v\in E_m.$$
    It is well known that $u\in E_m$ is a weak solution if and only if $u$ is a critical point of $J$. Since
$\overline{C_c^{\infty}}=H_ 0^m(\Omega)$, any weak solution of problem \eqref{aa}
    is also a solution in the distribution sense but it is not the case when $E_m=H_\vartheta^m(\Omega)$.
    Thanks to the assumption ($H$) and if
     $f\in C^{0,\alpha}_{loc}(\overline\O\times\mathbb{R})$, we can show that any weak solution of \eqref{a} (or \eqref{aa}) belongs
     to $C^{2m}(\overline\O)$
     (see \cite{VD}). This regularity allows us to obtain the lost part of Navier condition i.e.
    $\D^{j}u=0$ on $\partial \O$, for $\frac m2\leq j\leq m-1$ and to conclude that any weak solution is a
    classical solution of \eqref{a} or \eqref{aa}.

\medskip
   The application of critical points theory needs in general some compactness condition,
  known as the Palais-Smale condition (PS) which plays a central role :
  \begin{center}
If  $J(u_n)$ is bounded above and $J'(u_n)$  tends to $0$ in $E'_m$, then there exists a subsequence $u_{\phi(n)}$
  converging strongly in $E_m$.
  \end{center}
  In general such a condition requires the
  following standard assumptions \cite{AR}:
  \begin{itemize}
\item[(i)] Ambrosetti-Rabinowitz  condition:  There exist $\theta>2$ and $s_0>0$ such that
$$tf(x,s)\geq  \theta F(x,s) > 0\,\,\,\mbox{for} \,\,\, |s|>s_0\,\mbox{ and } x\in\O,.$$
\item[(ii)] Subcritical polynomial growth condition : There exist
\begin{align}
\label{psub} 1<p<\frac{N+2m}{N-2m} \;\;\mbox{if } N \geq 2m+1, \quad
\mbox{or  } 1 < p < \infty, \;\;\mbox{if } N \leq 2m,
\end{align}
$C>0$ and $s_0>0$ such that $|f(x,s)|\leq C|s|^{p}$ for $|s|>s_0$ and
$x\in\O$.
\end{itemize}
 (i) and (ii) are a standard Conditions in almost every work involving variational methods. However, (i) is so restrictive and it requires the following strong superlinear
  condition (SSL): There exist $C>0$ and $1<q\leq p $ such that $|f(x,s)|\geq
C|s|^{q}, C>0  $, for $x\in\O$ and $|s|$ large,(with $q=\theta-1$). When $m=1$, under the Dirichlet boundary condition, very few existence results have been established when $f$ satisfying (ii) and (i) is relaxed to (SSL)(see for example
\cite{DY, Y, WZ}). Nevertheless, (SSL) is also violated by many nonlinearities as for example $f(s)\sim as$ or $f(s)\sim as\ln(|s|)$ at infinity(where $a$ is a positive constant).
\subsection{Ambrosetti-Rabinowitz condition revised}
Let $u_n$ a (PS) sequence, condition (i) was made to get immediately $$|u_n|^2_m
    \leq  C_0(|u_n|_m +1),\,\, C_0>0.$$ . The aim here is to introduce a weaker condition than (i) which ensures that $u_n$ is bounded in $E_m$.
    The novelty in our approach is to write $J'(u_n)$ as a variational equation using the Riesz-Fréchet representation
    theorem as follows :
   there exists a unique $v_n\in E_m$ such that $$J'(u_n)\varphi=(v_n,\varphi)_m,\,\,\forall \varphi \in E_m\,\,\mbox{ hence }
   |J'(u_n)|_{E'_m}=|v_n|_m.$$
 Therefore $u_n-v_n$  could be seen as
 a weak solution in $E_m$ of
\begin{equation}\label{aaa}
  (u_{n}-v_{n},\varphi)_m= \int_{\O}h\varphi, \,\, \forall \varphi \in E_m,
   \end{equation}
where $h(x)=f(x,u_n(x))$. This device allows us to apply the well known $L^p$ elliptic regularity theory due to Agmon-Douglis-Nirenberg
\cite{ADN} :
\begin{lem}\label{l:1}
Let $p >1$, then there exists a positive constant $C_p$ such that for all $h\in L^p(\O)$, the following equation:
$$(u,\varphi)_m=\int_{\Omega}h\varphi,\,\,\, \forall  \varphi\in  E_m$$
admits a unique weak solution $u\in W^{2m,p}(\O)\cap E_m$.
Moreover, we have
$$ |u|_{ W^{2m,p}(\O)}\leq C_p|h|_{L^p(\O)}.$$
We denote by $L_p: L^p(\O)\longrightarrow W^{2m,p}(\O) $ the continuous linear operator defined by $L_p(h)=u$.
\end{lem}

We mention that our method is inspired by \cite{FLN, AH, HRS, KV, SO}, where a priori estimates involving the $L^{\infty}$ norm
was needed to establish existence results.

  \medskip
  Our first result is
\begin{pro}\label{pro:0}
 Let $f\in C (\overline\Omega\times \mathbb{R})$ satisfying ($H$) and \\
 ($H_1$) : There exist  $C>0$ and $s_{0}>0$
such that $$ C|f(x,s)|^{\frac{2N}{N+2m}}\leq sf(x,s)-
 2F(x,s) > 0,\,\, \forall\, |s|> s_{0},\,\,\,\ x\in \Omega.$$
 If  $u_n \in E_m$ satisfying
  $J(u_n)$  is bounded above and $J'(u_n)$ is bounded in $E'_m$, then  $u_n$ is bounded in $E_m$.
 \end{pro}
{\bf Remarks}\\
1) ($H$) and (i) imply ($H_1$). Indeed, $(i)$ implies $0<(1-\frac 2{\theta})sf(x,s)\leq sf(x,s) -2F(x,s)$,
for any $|s|>s_0$,  $x\in \O$. On the other hand $(H_1)$ implies  $C_0|f(s)|^{\frac{2N}{N+2m}}\leq sf(s)$ for all
$|s|>s_0$,  $x\in \O.$\\
2) Define $\xi(s)=\ln(\ln(\ldots \ln |s|)..),$ for $|s|$ large.
Then $f_\alpha(s)=s\xi^\alpha(|s|),$  satisfies  $(H_1)$ for all $\alpha > 0$.\\
 Of cours, $f(s)=as$ does not satisfy $H_1$, however
$f(s)=as-|s|^{\alpha-1}s$, $\frac{N-2m}{N+2m}\leq\alpha<1,$  $a>0$ and
 $f(s)=as -\frac{s}{\ln^{\alpha'}(|s|)},$ for $|s|$ large, $\alpha'>0$ and $a>0$ verify $H_1$. All these nonlinearities do not satisfy (i) and even (SSL).\\
3) If we assume the following very weak assumption $(H_0))$ : There exists $s'_0\geq s_0$ such that
$F(x,s'_0)$ and $F(x,-s'_0)$ are
positive  $\forall x\in \O$, then
 ($H_1$) and $(H_0)$ imply $(H)$. Thus, proposition \ref{pro:0} could be established only under $H_0$ and $H_1$.\\
4) The fact that $ L^r(\O)\hookrightarrow L^\frac{2N}{N+2m}(\O)\hookrightarrow E'_m$, for all
$ \frac{2N}{N+2m}\leq r\leq \infty$,
shows that the largest $L^p$ space including in $E'_m$ is $ L^\frac{2N}{N+2m}(\O)$. Consequently, if $u_n$ is a (PS) sequence, $(H_1)$
  imposes to
$f(x,u_n)$ to be bounded in $L^\frac{2N}{N+2m}(\O)$ which seems an optimal condition to
get the boundedness of $u_n$ in $E_m$.

\subsection{Subcritical polynomial growth condition revised}
After verifying that the (PS) sequence is bounded in $E_m$, we check the compactness of the nonlinear
  operator $K$ defined by (under assumption ($H$))
$$\begin{array}{lllll}
K:& E_m &\longrightarrow L^{\frac{2N}{N-2m}}(\O)) &\longrightarrow
 L^\frac{2N}{N+2m}(\O)& \longrightarrow E'_m\\
&  u &\longmapsto u &\longmapsto f(u) &\longmapsto f(u)
\end{array}$$
such that
$$K(u)v= \displaystyle\int_\Omega
    {f}(u)v, \quad \forall v\in E_m.$$
    Therefore (ii) and the Sobolev embedding ensure that $K$ is compact. However, (ii) is not satisfied
for the nonlinearity very close to the critical growth as for example
$$\dfrac {|s|^{\frac{4m}{N-2m  }}s}{(\ln(|s|+2))^q}, \quad \mbox{for } q>0.$$
We will use here the following strong subcritical condition weaker than (ii):
$$ \lim\limits_{t\to \infty}
\dfrac{f(s)}{|s|^{\frac{N+2m}{N-2m}}}=0,\,\,\,N\geq 2m+1. \leqno{(H_2)}$$
The condition $(H_2)$ seems to be optimal for ensuring that the operator $K$ is compact, since it is well known that
for the nonlinearity $f\sim |s|^{\frac{4m}{N-2m  }}s $ at $\infty$, $K$ is not
 compact.
\begin{pro}\label{pro:00}
 Let $f\in C (\overline\Omega\times \mathbb{R})$ satisfying ($H_2$). If $(u_n) \in E_m$ is bounded,
then there exists $u\in E_m$ such that $f(u_n)$ converges to $f(u)$ in $L^\frac{2N}{N+2m}(\O)$ up to a subsequence.
In other term, $K$ is compact.
\end{pro}
\subsection{Subcritical polynomial growth and Cerami condition}
 There is a variant of (PS) condition found by Cerami\cite{C}, noted by (C) :
 \begin{center}
 If $J(u_n)$ is bounded above and $(1+|u_{n}|_{m})|J'(u_{n})|_{E'_m}$ tends to $0$, then there exist a subsequence
 $u_{\phi(n)}$
  converging strongly in $E_m$.
  \end{center}
The condition (C) is stronger than (PS) condition, since we assume that $(1+|u_{n}|_{m})|J'(u_{n})|_{E'_m}$ converges to $0$ in addition. As
we have $2J(u_n)-J'(u_n)u_n\leq C_0$, which gives an advantage
that
the term $|u_n|_m$ does not appear
in the l.h.s., which is enough to begin the bootstrap process under the assumption $(H'_1)$ below.\\ Of course, to make this argument in work we need the subcritical polynomial growth
 condition (ii). More precisely, we see $J'(u_n)$ as a variational
 equation and we use the $L^p$ elliptic theory to obtain the uniform estimate for
$|u_n|_{L^{ p+1}}$ which yields that $u_n$ is bounded in $E_m$. we obtain then
 \begin{pro}\label{pro:000}
 Assume that $f$ satisfies $(ii)$ and
 $$\exists\; p_1>\sup(1,\frac{N(p-1)}{2mp}), \;C'_0, S_0 > 0 \mbox{ s.t. }C'_0|f(x,s|^{p_1}\leq sf(x,s)-2F(x,s),
 \,\,\forall\; |s|\geq s_0 \,\,x\in \O. \leqno{(H'_1)}$$
Then $J$ satisfies the Cerami condition.
\end{pro}
We remark that (ii) and $(H_1)$ imply $(H'_1)$ with $p_1=\frac{2N}{N+2m}$. Hence, if (ii) holds true,
($H'_1$) is weaker than $(H_1)$.\\
The following results are direct consequence of Propositions \ref{pro:0}-\ref{pro:000}.
\begin{pro}\label{p:0''}
 Let $f\in C (\overline{\Omega}\times \mathbb{R})$, Under assumptions $(H_1)$ and
$(H_2)$, $J$ satisfies the Palais-Smale condition.
  \end{pro}
\begin{pro}\label{pro:0000}
 Let $f\in C (\overline{\Omega}\times \mathbb{R})$ satisfying (ii) and
$(H'_1)$, then $J$ satisfies Cerami condition.
  \end{pro}
The assumptions $(H1)$ and $(H2)$ are complementary. In fact, $(H1)$ permits
to show that many nonlinearities close to the linear at infinity verify the (PS) condition; while (H2) is more interesting for
nonlinearities close to the critical growth. The combination of them would
permit to handle most nonlinearities with growth rate in the whole interval
$[1, \frac{N+2m}{N-2m}[$. Furthermore,  when $f$ is asymptotical linear at infinity, $H'_1)$ works better than $(H_1.$
\subsection{Ambrosetti-Rabinowitz Theorems revisited}
Thanks to the above propositions, we can improve the Ambrozetti-Rabinowitz theorems \cite{AR1, PR}.
 Denote by $\lambda_{1}=\lambda_{1,N}$ or $\lambda_{1}=\lambda_{1,D}$, the first eigenvalue of $(-\Delta)^m$ w.r.t.~the Navier
 and Dirichlet conditions.
 \begin{thm}\label{t:1}
  Let $f\in C (\overline\Omega\times \R)$ with $f(x, 0)=0$, $\forall x\in \O,$ assume $f$ satisfying ($H_1$),
($H_2$) and
$$\lim\limits_{t\to
\infty} \dfrac{f(x, t)}{t}  > \lambda_{1} \mbox{  uniformly in } \overline\Omega,\leqno{(H_3)}$$
$$\lim\limits_{t\to 0} \dfrac{f(x, t)}{t}  < \lambda_{1} \mbox{  uniformly in } \overline\Omega,\leqno{(H_4)}$$
Then $J$ has a nontrivial critical point obtained by the Mountain-Pass process.
 \end{thm}
 In \cite{FLN} de Figueiredo-Lions-Nussbaum considered a convex domain $\O$ and they proved existence results for
 the second order elliptic problem
  \begin{equation}\label{aaaa}
-\Delta u= f(u)\;\; in\; \Omega,\quad
 u=0\;\; on\;
\partial \Omega,
 \end{equation}
 under similar assumptions as Theorem\ref{t:1} ($m=1$) :
$$\lim\limits_{s\to \infty}
\dfrac{f(s)}{|s|^{\frac{N+2}{N-2}}}=0, \quad \lim\limits_{s\to \infty} \dfrac{f(s)}{s}>\l_1, \quad f(0)=0,\quad \lim\limits_{t\to
0} \dfrac{f(t)}{t}  < \lambda_{1} \leqno{(H_{FLN})}$$
with additional technical assumption. Similar conditions were made in \cite {ME, HRS}
where the domain $\O$ is a ball, to obtain radial solutions with prescribed number of zeros. Furthermore, it was
 conjectured in \cite{FLN} that for general domain, if $f$ is locally Lipschitizian
satisfying {\it only} $(H_{FLN})$, the equation \eqref{aaaa} has a classical nontrivial positive solution.
When $f$ does not depend on $x$, they partially proved this conjecture in Theorem 2-3 in \cite{FLN}, under some global restrictive
assumption $(26)$.\\
 We point out that Theorem \ref{t:1} here gives a partial answer for this conjecture under the assumption $(H_1)$ (with $m=1$), also
  for the poly-harmonic equation \eqref{a} when $f$ is a nonnegative function and $m\geq 2.$. Moreover, our solutions have
  a min-max structure.
 \begin{cor}\label{c:1}
 Let $m=1$ and $f\in C (\overline\Omega\times \R)$ with $f(x, 0)=0$, $\forall x\in \O.$ Assume that $f$ satisfies
 ($H_1$)-($H_4$), then \eqref{aaaa} has a positive solution.
  \end{cor}
 It is well known that if $f$ is odd then (i) and (ii) allow removal of any condition near $0$ for $f$ to obtain infinitely solution of \eqref{a} or \eqref{aa}. However, under $(H_1)$ and $(H_2)$ we need to add the following assumption at $0$ $(H'_4)$ :$\lim\limits_{t\to
0} \dfrac{f(t)}{t}  <\infty$. We get
  \begin{thm}\label{t:1'}
 Let $f\in C (\overline \Omega\times \mathbb{R})$ be an odd function. Assume that $f$ satisfies ($H_1$)-($H_3$) and $(H'_4)$ then
 $J$  admits infinitely many distinct pairs $(u_n,-u_n)$, $ n\in \mathbb{N}$ critical points of $J$,
 obtained by the symmetry Mountain Pass theorem. Moreover, $J(u_n)$ is unbounded.
  \end{thm}
  Consider the following examples
$$f(s)\dfrac {|s|^{\frac{4m}{N-2m  }}s}{(\ln(|s|+2))^q}, \quad f(s) = s\ln\left(\ln..\ln(|s| +C)..\right) \mbox{  for $C$ large}.$$
This  nonlinearities satisfy ($H_1$)-($H_3$ and $(H'_4)$\\
{\bf Remarks}\\
1)Under $H_1$-$H_3$ and $(H'_4)$ we can use the linking theorem to obtain a non trivial critical point of $J$.\\
2)  In the case of Hamiltonian systems, existence results could be established if we relax the Ambrozetti-Rabinowitz condition
by the following condition
$H''_1$: $ |H(z)|\leq H_z(z).z-2H(z)$ for $|z|$ large.\\
3) All these existence results could be established if we replace $H_1$ by $H'_1$ and $H_2$ by $(ii)$
and we observe that $f(s)=as-|s|^{\alpha-1}s$,
$\alpha<1$ and $a>\lambda_1$  satisfies $H'_1$, $(ii)$, $H_3$ and $H_4$. Hence, the mountain pass theorem applied.
4)In\cite{CM} and when $m=1$,
the authors assume the following:\\
There exist $r\in[1,\frac{N+2}{N-2}[$ and $C_0>0$, such that $f(x,s)\leq C_0(|s|^p +1)$, $\forall (x,s)\in\O\times\mathbb{R}$.
In addition, there exist $\mu>r-1$, if $N=1, 2$ and $\mu>{N(p-1)}2$, if $N\geq 3$ and two positive constants $C'_0$, $S_0$
such that $$C'_0|s|^{\mu}\leq sf(s)-2F(s), \,\,\,\forall s\geq s_0..$$ Many interesting existence results were proved involving
the asymptotical linear case or resonance case at $\infty$( see also \cite{CM2 ZZ}). However,
only one superlinear case was studied under an additionnel global assumption
namely: $F_4$ $f(x,s)\geq\frac{\l_{k-1}s^2}2$ {\it for all} $ (x,s)\in \O\times\mathbb{R}$.\\
\section{\bf Weaken subcritical case}
To motivate our next result we consider the following instructive example.\\
Set $f(s) = \gamma s^q+s^p(1+\sin(\ln(\ln(s)))$, $s\geq 0,$ and $f(s)=0$, $s<0$ where $1<q\leq  p$ satisfy \eqref{psub} and $\gamma\geq 0$
 .\\
-If $q>p$, then $f$ satisfies condition (i).\\
-If $q=p$ and $\gamma(p+1)>2(\gamma+1)$, then $f$ satisfies condition (i).\\
-If $q=p$ and $\gamma> 1$, then $f$ satisfies the strong superlinear condition (SSL)\\
- If $\gamma=0$, then $f$ does not satisfy even the weak superlinear
condition $H_3$ since $f(l_n)=0$ where
$l_n=\exp\exp((2n-\frac12)\pi)$.  Hence, $\dfrac{f(t)}{t}$ and $f$  {\bf do not have any limit at infinity}, this case will be called the weaken subcritical case. \\
 In the following, we improve Theorem (1.1) obtained in \cite{AH}. In fact, we will remove assumption $(f_2)$ to prove theorem (1.1) for all $m\geq 1$
 recall that Theorem (1.1) was established without assuming the weak superlinear condition $H_3$ or any other assumption on the behavior of $ \dfrac{f(t)}{t}$ at infinity as it is the case in Ambrosetti-Prodi-type problems( see \cite{K}, for example). let $ \l_1 $ be the first eigenvalue of $(-\Delta^)m$ under the Navier
boundary condition, $\R_+ = [0, \infty)$, consider the problem \eqref{a} with $f(x,s)=f(s)$.
 \begin{thm}\label{t:1''}
Problem \eqref{a} has at least one positive solution $u\in C^{2m}(\O)\cap C^{2m-1}(\ov\O)$ provided that $f$ satisfies\\
$(f_1)$ $f \in C^1 (\R_+, \R_+)$, $f(0)=0$, $ f'(0) <
\lambda_1.$ Moreover, there exist $C>0$, $T>0$ and $ p$ satisfying \eqref{psub} such that $|f'(s)|\leq C(|s|^{p-1}+1)$. \\
$(f_2)$ there exist $s_n \to \infty$ and $\mu>0$ such that
$f(s_n)\geq \mu
s_n^p.$\\
$(f_3)$ $\,\,\,\,\,\,\,\,\,\,\,\,\,\, \lim\limits_{t\to
+\infty}\dfrac{sf'(s)-pf(s)}{s^p}=0.$
 \end{thm}
 We remark that the example above corresponding to $\gamma=0$ :  $f(s)=s^p(1+\sin(\ln(\ln(s)))$, $s\geq 0,$ satisfies $(f_1)$-$(f_4).$\\
  This paper is organized as follows. In section 2, we give the proof of propositions \ref{pro:0}, \ref{pro:00}, \ref{pro:000}
 . Section 3 is devoted to the
proof of Theorem \ref{th:1}, \ref{th:1'}, \ref{th:1''} and corollary`\ref{c:1}. In section 4, we give the proof of
Theorem \ref{l:2}, \ref{l:3} and proposition\ref{p:1}.\\
\section{ \bf proof of propositions \ref{pro:0}, \ref{pro:00}, \ref{pro:000} }
 { \bf proof of propositions \ref{pro:0}}\\
Let $u_n$ be a sequence such that $|J'(u_n)|_{E'_m}$ is bounded and $J(u_n)$ is bounded above.  The Riesz-Fr\'echet representation theorem ensures the existence
of $v_n\in E_m$ such that $J'(u_n)\varphi=(v_n,\varphi)_m,$  $\forall\;\varphi \in E_m$ with $|J'(u_n)|_{E'_m}=|v_n|_m$. Consequently, we have $$(u_n-v_n,\varphi)_m=\int_{\O}f(u_n)\varphi,\,\, \forall \varphi\in E_m.$$  Moreover, we have $\int{\O}(f(x,u_n)u_n -2F(x,u_n))=2J(u_n)- J'(u_n)u_n\leq C_0(|u_n|_m+1), $ then $(H_1)$ implies that $|f(x,u_n)|_{L^{\frac{2N}{N+2m}}(\O)}\leq C'_0(|u_n|_m+1). $ Thus, by lemma\ref{l:1} we deduce that $$|u_n-v_n|_{ W^{2m,\frac{2N}{N+2m}}(\O)}\leq C|f(x,u_n)|_{L^{\frac{2N}{N+2m}}(\O)} \leq C''_0(|u_n|_m+1). $$  Taking in account that  $W^{2m,\frac{2N}{N+2m}}(\O)\hookrightarrow W^{m, 2}(\O)=H^m(\O),$ then we derive that $$|u_n-v_n|_m^{\frac{2N}{N+2m}}\leq C_0(|u_n|_m +1).$$ Finally, since $|v_n|_m=|J'(u_n)|_{E'_m}$
is bounded and $\frac{2N}{N+2m}>1$, we get $|u_n|_m$ is bounded.\\
{ \bf proof of propositions \ref{pro:00}}\\
By $H_2$ we deduce that for every $\epsilon>0$ there exists a positive constant $c_{\epsilon}$ such that
\begin{equation}
\label{e:11}
|f(x,s)|^{\frac{2N}{N+2m}}\leq \epsilon |s|^{\frac{2N}{N-2m}}+c_{\epsilon},\,\, \forall (x,s)\in \overline\O\times \R.
\end{equation}
 Since $u_n$ is bounded in $E_m$, and $f$ is a continuous function in $\overline{\Omega}\times \mathbb{R}$,  then there exists a subsequence $u_{\psi(n)}$,  $u\in E_m$ and a positive constant $C_0$ such that $$f(x,u_{\psi(n)}) \rightarrow f(x,u) \mbox{  a.e in } \O,\,\,\,  |f(x,u)|<\infty \mbox{ a.e in  } \O$$ and
 \begin{equation}
\label{e:11'}
 2^{\frac{2N}{N+2m}}\int_{\O}|u_{\psi(n)}|^{\frac{2N}{N-2m}}+|u|^{\frac{2N}{N-2m}}\leq C_0.
 \end{equation}
  Therefore, Egorov theorem implies that there exists $\O_{\epsilon}\subset \O$ such that mes($\O\setminus \O_{\epsilon})<\frac{\epsilon}{c_{\epsilon}}$ and $|f(u_{\psi(n)})-f(u)|_{L^{\infty}(\O_{\epsilon})}$ converges to $0,$ when $n$ tends to $+\infty,$ which also implies that $$\int_{\O_{\epsilon} }|f(u_{\psi(n)})-f(u)|^{\frac{2N}{N+2m}},$$ converges to $0,$ when $n$ tends to $+\infty$ .
  On one hand, we have $$\int_{\O\setminus\O_{\epsilon}}|f(u_{\psi(n)})-f(u)|^{\frac{2N}{N+2m}}\leq 2^{\frac{2N}{N+2m}}(\int_{\O\setminus \O_{\epsilon}}(|f(u_{\psi(n)})|^{\frac{2N}{N+2m}}+|f(u)|^{\frac{2N}{N+2m}}),$$ then \eqref{e:11}, \eqref{e:11'} and the fact that  mes($\O\setminus \O_{\epsilon})<\frac{\epsilon}{c_{\epsilon}}$ imply that
  \begin{equation}
\label{e:12}
  \int_{\O\setminus\O_{\epsilon}}|f(u_{\psi(n)})-f(u)|^{\frac{2N}{N+2m}} \leq (C_0 +2)\epsilon.
  \end{equation}
On the other hand, since $$\int_{\O_{\epsilon} }|f(u_{\psi(n)})-f(u)|^{\frac{2N}{N+2m}},$$ converges to $0,$ we derive that
 there exists $N_{\epsilon}$ such that for $n>N_{\epsilon}$ we have  $$\int_{\O_{\epsilon} }|f(u_{\psi(n)})-f(u)|^{\frac{2N}{N+2m}}\leq  \epsilon.$$
 To conclude,  thanks to \eqref{e:12} and the least inequality , we obtain for $n>N_{\epsilon}$,  $$\int_{\O}|f(u_{p(n)})-f(u)|^{\frac{2N}{N+2m}}\leq(C_0+3)\epsilon.$$ Hence, the result follows.\\
{ \bf proof of propositions \ref{pro:000}}\\
Let $u_n$  a (C) sequence, thanks to (ii), if we prove that $u_n$ is bounded in $L^{p+1}(\O),$ then $u_n$ is bounded in $E_m.$   Applying again our approach as in the proof of proposition \ref{pro:0} : there exists $v_n\in E_m$ such that $$(u_n-v_n,\varphi)_m=\int_{\O}f(u_n)\varphi,\,\, \forall \varphi\in E_m,$$ and $|J'(u_n)|_{E'_m}=|v_n|_m$. Moreover, since $|J'(u_n)|_{E'_m}|u_n|_m$ is bounded and $J(u_n)$ is bounded above, then there exists a positive constant $C_0$ such that $$\int_{\O}(f(x,u_n)u_n -2F(x,u_n))=2J(u_n)- J'(u_n)u_n\leq C_0 .$$  As a consequence, $u_n$ is bounded in $L^{p+1}(\O)$ if and only if $u_n-v_n$ is bounded in $L^{p+1}(\O).$ Furthermore,  $(H'_1)$ implies that $f(u_n)$ is bounded in $L^{p_1}(\O)$, using lemma\ref{l:1} to obtain $(u_n-v_n)$ is bounded in $W^{2m,p_1}(\O)$, by Sobolev embedding we have $(u_n-v_n)$ is bounded $L^{p^{*}_1}(\O)$ where $$p_1^{*}=\frac {Np_1}{N-2mp_1},\mbox{ if } N>2mp_1 \mbox{ or }p_1^{*} =+\infty \mbox{ if }  N<2mp_1.$$
. Therefore, if $p+1\leq p_1^*$, we are done. If not, observe that the condition $p_1>\frac{N(p-1)}{2mp}$ implies that $\frac{p_1^{*}}{p}>p_1$, and by $(ii)$, we get $f(u_n)$ is bounded in $L^{\frac{p_1^{*}}{p}}(\O)$. By consequence, we iterate the $L^p$ regularity, to obtain $p_{k+1}=\frac{p_{k}^*}p$ and $u_n\in ^L{p_{k+1}^*}(\O)$, where $p_{k+1}^*=\frac {Np_{k+1}}{N-2mp_{k+1}}$ if $ N>2mp_{k+1}$ and $p_{k+1}^{*} =+\infty$ if $ N<2mp_{k+1}.$ We argue by contradiction, we suppose that $$p_{k+1}^*=\frac {Np_{k+1}}{N-2mp_{k+1}}<p+1,\,\, \forall\,\, k\in \mathbb{N}$$.\\
 {\bf Case 1 $ p>1$}. Set $r_k=\frac 1{p^{*}_k}$, $k\in\mathbb{N}\setminus\{0\}$ and $r_0=\frac 1{p_1}$, we have, $r_{k+1}=pr_k-\frac{2m}{N}$, then $$r_{k+1}=p^k(\frac 1{p_1}-\frac{2mp)}{ N(p-1})+\frac{(2m)}{ N(p-1)}$$. Since $p_1>\frac{N(p-1)}{2mp}$ and $p>1$ we derive that $r_k$ tends to $-\infty$ but we have $r_k>\frac {1}{p+1},\,\, \forall k\in \mathbb{N} $, a contradiction.
To conclude, there exists $k_0\mathbb{N}\setminus\{0\}$ such that $u_n$ is bounded  in $L^{p+1}[\O).$\\
 {\bf Case 2 $ p=1$}. We have $p_{k+1}^*=\frac {Np_1}{N-2kmp_1}$, hence for $k$ large enough we get $p_{k+1}^*=+\infty$, then we deduce that $u_n$ is bounded in $L^{2}[\O).$\\
Finally, the fact that $u_n$ is bounded in $E_m$ and (ii) imply that $u_n$ has a subsequence converging strongly in $E_m$.\\
\section{ \bf proof of theorem \ref{t:1}, \ref{t:1'}, \ref{t:1''} and corollary \ref {c:1}}
{ \bf proof of theorem \ref{t:1}}\\
Assumption $(H_2)$ and $(H_4)$ imply that there exist $r>0$ and $\alpha>0$ such that $J(u)\geq \alpha$, for all $u\in E_m, \,\, |u|_m=r$. $(H_1)$ and $(H_2)$ imply the $(PS)$ condition. Therefore, to prove theorem \ref {th:1'}, we need only to verify that there exists $v_0\in E_m$ such that $J(v_0)<0$. We have $H'_3$ implies $\lim\limits_{s\to \infty} \dfrac{F(s)}{s^2} >\frac{\l_1}2.$ Thus, there exist $ \epsilon_0>0$ small enough and a positive constant $c_0$ such that $ |F(s)| \geq  (\frac{\l_1}2+\epsilon_0) |s|^ 2-c_{ 0}$. Let $\varphi_1$ be the eigenfunction associated to $\l_1$, set $v_0=\beta \varphi_1$ , $\beta>0$  . We derive that $J(v_0)\leq \frac{\beta^2}2(|\varphi_1|^2_m -\l_1\int_{\O} \varphi_1^2)-\beta^2\epsilon_0\int_{\O} \varphi_1^2+C_0$mes$(\O)\leq -\beta^2\epsilon_0\int_{\O} \varphi_1^2+C_0$mes$(\O)$,  . Hence, for $\beta$ large enough, we get $J(v_0)<0$.\\
 {\bf proof of corollary \ref{c:1}}\\
 When $m=1$, In order to obtain a positive weak solution, we may truncate $f$ above $0$ by denoting $f^{+}(x,s)=0, \,\,\,s<0$ and $f^{+}(x,s)=f(x,s), \,\,s\geq 0$. If $f$ verifies $H_1$-$H_4$ for $s>o$, the corresponding energy functional $J$ satisfies the hypothesis of mountain pass theorem. Since $u^+=\sup(u,0)\in H^1_0(\O), \;\; u^-=\sup(-u,0)\in H^1_0(\O) $, $\nabla u^=\nabla u^{-}+\nabla u^{+}$ e.a in $\O$ and $\int_{\O}\nabla u^{+}.\nabla u^{-}=0$ \cite{OK}, if we multiply equation \eqref{aaaa} by $u^-$ and we integrate by part, we derive that
   $\int_{\O}|\nabla u^{-}|^2=0$, then $u$ is a positive weak solution. Finally, The $L^p$-estimates employed by Brezis and Kato \cite{BK} implies that $u$ is classical solution.  A very interesting {\it open question} is how can us remove assumption $H_1$ ?. We point out that positive solution for higher order equations is harder to exhibit; loosely speaking, this is due to fact that the decomposition $u= u^++u^-$ is no longer available in $E_m$ for $m\geq 2$.\\
 {\bf proof of theorem \ref{t:1'}}\\
Denote by $0<\lambda_1<\lambda_2\leq \lambda_3 . . \leq \lambda_k$ the eigenvalues of $(-\Delta)^m$ with Dirichlet or Navier boundary condition and let be $\varphi_k$ the corresponding eigenfunctions. We claim that for $k_0$ sufficiently large there exist $r>0$ and $\alpha>0$ such that for all $u\in E^{+}_m :=\mbox{span}\{\varphi_k; k\geq k_0\}$ with $|u|_m\geq r$ there holds $J(u)\geq \alpha$. Indeed, $(H_2)$ and $(H'_4)$ imply that
$$|F(x,s)| \leq  |s|^{\frac{2N}{N-2m}}+C|s|^2,\,\, \forall (x,s)\in \overline\O\times \R.$$
Hence, $J(u)\geq \frac12|u|_m^2 - \int_{\O} |u|^{\frac{2N}{N-2m}}-C\int_{\O}|u|^2.$ Taking in account that $u\in E^+_m$, we get
$$J(u)\geq (\frac12 -C_0|u|_m^\frac{4m}{N-2m}-C_0\lambda_{k_0}^{-1})|u|_m^2,$$ (where $C_0>0$ is independent of $u$). Choose $k_0$ large enough such that $\frac12-C_0\lambda_{k_0}^{-1}=\frac14$, then we get $J(u)\geq (\frac14 -C_0|u|_m^\frac{4m}{N-2m})|u|_m^2$ hence we may take $r=\frac1{C_0^\frac{N-2m}{8m}}$ which implies that $\alpha=\frac18r^2.$ Thus, theorem \ref{t:1'} is well proved.\\
 {\bf Proof of Theorem \ref{th:1''}}
We consider the classical truncated function:
$$f_n(s)=\left\{\begin{array}{ll}
0  & \hbox{if } s\leq 0 \\
f(s) & \hbox{if }0\leq s\leq s_n \\
 f(s_n)-\frac1p s_nf'(s_n)+\frac1p\displaystyle\frac{f'(s_n)}{s_n^{p-1}}s^p & \hbox{if } s\geq
 s_n.
\end{array}
\right.$$ We have $f_n\in C^0(\R, \R_+)$ and $f_n\in C^1(\R_+,
\R_+)$, We can check easily that $(f_1)$, $(f_3)$ and the definition
of $f_n$  imply that
\begin{align}
\label{f1,2unif}\mbox{ there exists}\;\; C'>0 \;\; \mbox{such that}
  |f'_n(s)|\leq C'(1+s^{p-1}),\quad   0\leq f_n(s)\leq C'(1+s^p ),\quad
\forall\, (s,n)\in \R_+\times\N  .
\end{align}
Furthermore, it is not difficult to see that by $(f_3)$ and the
definition of $f_n$, we have

\begin{align}
\label{f4unif}
\frac{sf_n'(s) - pf_n(s)}{s^p} \;\; \mbox{tends uniformly to $0$, as } s \to +\infty.\end{align}

\medskip
Next, we consider the truncated problem
\begin{equation}\label{e:i}
\left\{\begin{array}{lll} \Delta^2 u=f_n(u)\ \ \  \mbox{in}\,\, \O,\\
\,\ \ u=\D u=0\ \ \ \mbox{on}\ \ \partial \O.
\end{array}
\right.
\end{equation}
The mountain pass-theorem shows that equation \eqref{e:i} has a
positive classical solution $u_n$ (see \cite{VD}). We claim that
$\|u_{n_0}\|_{L^{\infty}(\O)}\leq s_{n_0}$, for some $n_0$. Then
$u_{n_0}$ is a solution of
 \eqref{e:a} and Theorem \ref{th:1} follows.

\medskip
 We argue by contradiction and assume that $ \|u_n\|_{L^{\infty}(\O)} >s_n$, for all
 $n$.
Then $\|u_n\|_{L^{\infty}(\O)} \to \infty$ as $n$ tends to $\infty$.
  Define
$\b_1=\frac{2m}{p-1}$,

 $\l_n^{-\b_1} =  \|u_n\|_{L^{\infty}(\O)} $, and $x_n\in\O$ such that  $u_n(x_n) = \l_n^{-\b_1}$,

$\tilde u_n(y)=\l_n^{\beta_1} u_n(x_n+\l_ny)$, $ y\in\widetilde \O_n = \l_n^{-1}(\O - x_n)$

. we have $\l_n \to 0$  and  $\|\tilde u_n\|_{L^{\infty}(\widetilde\O_n)}=\tilde u_n(0)=1$.
  Moreover, $ \tilde u_n$ satisfies
$$(- \Delta)^m \tilde u_n = \l_n^{p\beta_1}f_n(\l_n^{-\beta_1}\tilde u_n)\,\,\, \mbox{in}\,\,\,\widetilde\O_n.$$
Set $g_n(y): = \l_n^{p\beta_1}f_n(\l_n^{-\beta_1}\tilde u_n(y))$.
Then,  \eqref{f1,2unif} implies that
$\|g_n\|_{L^{\infty}(\widetilde\O_n)}$ is bounded.
\medskip
{\bf Case $1$.} Up to a subsequence, there holds $\l_n^{-1}{\rm
dist}(x_n,\partial\O)\to \infty$.
\smallskip
Then, the domain $\widetilde \O_n = \l_n^{-1}(\O - x_n)$ converges to
the entire space $\R^N$. Then, for all $R>0 \,\, \tilde u_n$ is
defined in $B_{3R}(0)$ for $n$ large enough. Applying Corollary $6$
in\cite{RW} (page 809), we obtain $\|\tilde u_n\|_{W^{4,q}(B_{2R})}$
is bounded for all $q>1$. Hence, $\tilde u_n$ is bounded in $C^{2m-1,\a}(B_{2R})$  and the mean value theorem together with
\eqref {f1,2unif} imply that $g_n$ is bounded in $C^{0,\a}(B_R)$ for any $R > 0$
  . Thus, using the standard Schauder estimates we derive that $\tilde u_n$ is bounded in $C^{2m,\a}(B_R)$.
  By the diagonal process, up to a subsequence,  the Ascoli-Arzela's
  theorem implies that
  $\tilde u_n$ and  $\tilde v_n$ converge to
 $u$
 in $C^{2m, \a'}_{loc}(\R^N), 0<\a'<\a$,  and $g_n$
converges to $g$ in $C^{0, \a'}_{loc}(\R^N)$. Moreover, $u, g
\geq 0$  in $\R^N$, $u(0) = 1$, and
$$  (-\Delta u)^m = g \quad \mbox{in }\R^N.$$
As $u \not\equiv 0$, the strong maximum principle ensures that $u
> 0$ in $\R^N$. Hence,  $\l_n^{-\beta_1}\tilde u_n$
tends to infinity uniformly on each compact set of $\R^N$. Set $Q_n
= g_n\tilde u_n^{-p}$ and $Q = gu^{-p}$, we have $Q_n$ converges to $Q$ in $C^0_{loc}(\mathbb{R}^N$
\medskip
 . Furthermore,
\begin{align*}
|\nabla Q_n| = \left|\frac {\tilde u_n\nabla g_n - pg_n \nabla\tilde  u_n}{\tilde u^{2p}_n}\right| =
 \left|\frac {f_n'(\l_n^{-\beta_1}\tilde u_n) - pf_n(\l_n^{-\beta_1}\tilde u_n)}{(\l_n^{-\beta_1}\tilde u_n)^{p - 1}}\right|\times
  \frac{|\nabla\tilde u_n|}{\tilde u_n^p}
\end{align*}
Since $u>0$ then for any $R>0$, there exists $C_R>0$ such that $ \tilde
u_n(x)>C_R, \forall x\in B_R$ uniformly in $n$,  which implies that  $\l_n^{-\beta_1}\tilde
u_n$ tends to $\infty$ uniformly on $B_R$ and
$\frac{|\nabla\tilde u_n|}{\tilde u_n^p}$ is bounded in $C^0(B_R)$.
 Thus, by
\eqref{f4unif} we derive that
$$\left|\frac {f_n'(\l_n^{-\beta_1}\tilde u_n) - pf_n(\l_n^{-\beta_1}\tilde u_n)}{(\l_n^{-\beta_1}\tilde u_n)^{p -
1}}\right|,$$ tends to zero in $C^0 (B_R)$. Hence, we can conclude that $\nabla Q_n \to
0$ in  $C^0 (B_R), \forall R>0 $. Moreover, since $ Q_n$ could be
written as $ Q_n(x)=\int_0^{x_i}\frac{ \partial  Q_n}{\partial
x_i}(x_1,..  t, x_{i+1},.. ,x_N)dt- Q_n((x_1,..  0, x_{i+1},..
,x_N))$, $\forall i\in\{1,2....,N\}$ then $ Q_n(x)$ converges to
$Q((x_1,.. 0, x_{i+1},.. ,x_N))$. By iteration, we derive that $
Q_n(x)$ converges to $Q(0)$. Hence $Q(x)=Q(o)\geq 0$.  More precisely, by $(f_2)$ and $(f_3)$, we deduce that $Q(0)>0$. Indeed, recall that
$\l_n^{-\b_1}= \|u_n\|_{L^{\infty}(\O)} >s_n$ and $\tilde u_n(0)=1$,
then
   $$Q_n(0)=\l_n^{p\beta_1}f_n(\l_n^{-\beta_1})=\l_n^{p\beta_1}(f(s_n)-\frac1p s_nf'(s_n))+\frac1p\displaystyle\frac{f'(s_n)}{s_n^{p-1}}.$$
Thus, $(f_3)$ implies that $\l_n^{p\beta_1}(f(s_n)-\frac1p s_nf'(s_n))$ tend to $0$, also by $(f_2)$ and $(f_3)$ we derive
  that there exist  $0<\mu'<\mu$ and $N_0\in \mathbb{N}$ such that $f'(s_n)\geq \mu' s_n^{p-1}, \forall n>N_0.$  Hence, we have
 $$ 0< \frac{\mu' }{p}\leq \lim \frac1p\dfrac{f'(s_n)}{s_n^{p-1}}$$
 and
$$\lim Q_n(0)=Q(0)  >\frac{\mu' }{p}>0$$. So, we obtain a positive
solution $u$ of $(-\Delta)^m u = Q(0)u^p$ in
  $\R^N$, where $Q(0)$ is a positive constant and $p$ is subcritical, which is impossible (see
  \cite{WX}).
\medskip
 {\bf Case  2.} $\l_n^{-1}{\rm dist}(x_n,\partial\O)$ is bounded. Up
to a subsequence, assume that $\l_n^{-1}{\rm dist}(x_n,\partial\O)
\to d \geq 0$. Since $u_n + v_n = 0$ on $\partial \O$,  a standard
elliptic estimate prove that $d > 0$, (see, for example \cite{RW}).
By flattening the boundary through a local change of coordinate we
may assume that near $x=\lim x_n$ the boundary is contained in $x_N
=0.$ Now, using standard scaling and translation argument, as for
the above case and applying Corollary $6$ in\cite{RW} (page 809)
(with Navier boundary condition) to obtain a nontrivial classical
solution for $(-\Delta(^m u = Qu^p, u(0)=1 $
   in the half space $X_N>-d$ with $u = \Delta u =. . =(\Delta u)^{m-1} 0$ on $x_N=-d$. Finally, we can prove that $Q$
   is a positive constant, which is again impossible thanks to Sirakov's result \cite {KV}). So, we are done.\\


\begin{thebibliography}{777}
\bibitem{ADN}
S. Agmon, A. Douglis, L. Nirenberg, \emph{ Estimates near the boundary for solutions of elliptic  partial
differential equations, satisfying general boundary conditions} Comm. Pures Appl. Math \textbf{12} 1959

\bibitem{AR}
A. Ambrosetti and P. H. Rabinowitz,\emph{Dual variational methods in critical point theory and applications,}
  J. Funct. Anal.  \textbf{14} (1973), 349

 \bibitem{BK}
H. Brezis and T. Kato, \emph{Remarks on the schrodinger operator with singular compact potentials,} J. Math. Pures Appl.
 \textbf{58} (1978), 137-151.

\bibitem{C}
G. Cerami, \emph{An existence criterion for the critical points in unbounded manifolds,} J. Rend. Sci. Mat. Fis. Geol. 112 (1978), 332-336

\bibitem{CM}
C.G. Costa and C.A. Magalhaes,\emph{ A unified approach to a class of Strongly indefinite functionnals ,} J. Differntial Equations \textbf{125} (1996)
\bibitem{CM2}
C.G. Costa and C.A. Magalhaes,\emph{ A variationnel approach to subquadratic perturbation of elliptic system a class ,} J. Differntial
Equations \textbf{111} (1994)

\bibitem{ME}
M.J. Esteban, \emph{Multiple solutions of semilinear elliptic problems in a ball ,} J. Diff. Equ. \textbf{57} (1985), 112-137.

\bibitem{FLN}
D.G. de Figueiredo, P. L. Lions and R. Nussbaum, \emph{A priori
estimates and existence of positive solutions of semilinear elliptic
equations,} J. Math. Pures Appl. \textbf{61} (1982), 41-63.

\bibitem{DY}
D.G. de Figueiredo and J. Yang, \emph{On a semilinear elliptic problem without (PS) condition,} J. Differntial Equations \textbf{187} (2003), 412-428.

\bibitem{GGS}
F. Gazzola, H. Christoph Grunau and G. Sweers,  \emph{Polyharmonic boundary value problems} Springer.

\bibitem{AH}
  A. Harrabi, \emph{ On a fourth-order superlinear elliptic problem without (PS) condition,}  accepted in Adanced nonlinear studies (2013).

  \bibitem{HRS}
A. Harrabi, S. Rebhi and A. Selmi, \emph{Existence of radial solutions with prescribed number of zeros of elliptic solutions and their Morse index,} J. Differntial Equations \textbf{251}(2011),  2409-2430.
\bibitem{OK}
O. Kavian, \emph{Intoduction à la théorie des points critiques,} Mathématiques et Applications \textbf{13} (1991).

\bibitem{K}
E.Koizumi, K Schmitt, \emph{Ambrosetti-Prodi-type problems for quasilinear elliptic problems} Diff. Int. Equa.\textbf{18:3} (2005)  241-262.
\bibitem{PR}


\bibitem{RW} W. Reichel and T. Weth, \emph{A priori bounds and a Liouville theorem on a half space for higher-order elliptic Dirichlet problems,} Math. Z. {\bf 261} (2009), 805-827.
\bibitem{KV}

\bibitem{K}
E.Koizumi, K Schmitt, \emph{Ambrosetti-Prodi-type problems for quasilinear elliptic problems,} Diff. Int. Equa.\textbf{18:3} (2005)  241-262.
\bibitem{PR}
P.H. Rabinowitz, \emph{Min-Max methods in critical point theory with applicatios to differential equations,} Conf. Board of the Mathematical sciences.\textbf{53}(1984).
\bibitem{RR}
M. Ramos and P. Rodrigues \emph{On a fourth order superlinear
elliptic problem,} Electron. J. Diff. Eqns, Conf. \textbf{06}
(2001), 243-255.
114-133.
\bibitem{MS}
M. Struwe, \emph{ Variatinal methods, Applications to nonlinear partial differential equations and Hamiltonian systems,} Fourth edition, volume (34) Springer.

\bibitem{RW} W. Reichel and T. Weth, \emph{A priori bounds and a Liouville theorem on a half space for higher-order elliptic Dirichlet problems,} Math. Z. {\bf 261} (2009), 805-827.
\bibitem{KV}
B. Sirakov, \emph{Existence results and a priori bounds for higher order elliptic equations and systems,} J.math.Pures.Appl \textbf{9} 89 (2008), no. 2, 114-133.
\bibitem{SZ}
J. Su and L. Zhao \emph{ An elliptic resonance problem with multiple solutions,} J.math.Anal.Appl.
\textbf{319}  (2006), pp 604-616.

\bibitem{SO}
R. Soranzo, \emph{A priori estimates and existence of positive
solutions of superlinear polyharmonic equation,} Dynamic System and
Appl. \textbf{3} (1994), 465-488.

\bibitem{VD}
 R. Van Der Vorst, \emph{ The best constant for the embedding $L^{\frac{2N}{N-4}}$ in $H^2(\O)\cap H_0^1(\O)$,}  Journal of Diff and Integ Equat \textbf{6} (1993), no. 2. 259-273.
\bibitem{WX}
J. Wei X. Xu \emph{A classification of solution of higher order conformally
invariant  equations,} Math. Ann
 \textbf{313}(1999), 207-228.

\bibitem{Y}
 J. Yang, \emph{Existence of solutions of semilinear elliptic
problem without (PS) condition,} Proc. Amer. Math. Soc. 132(5) (2004), 1355-1366. MR2053340 (2005b:35048).

\bibitem{Y2}
 J. Yang , \emph{Corrigendum to Existence of solutions of semilinear elliptic
problem without (PS) condition,} Proc. Amer. Math. Soc. S
0002-9939(2011)10909-2 Article electronically published on June 2,
2011
\bibitem{WZ}
 W. Zou, \emph{Variant fountain theorems and their applications,} Manuscripta Math. 104 (2001),
343 - 358. MR 2002c:35081.

\end{thebibliography}
 \end{document}